\documentclass[12pt]{amsart}

\usepackage{url}
\usepackage{amssymb}
\usepackage{nicefrac}
\usepackage{tikz}
\usetikzlibrary{shapes,positioning,intersections,quotes}
\usepackage[utf8]{inputenc}
\usepackage{float}

\usepackage{enumitem}

\usepackage{graphicx}

\makeatletter
\@namedef{subjclassname@2020}{%
  \textup{2020} Mathematics Subject Classification}
\makeatother

\usepackage[T1]{fontenc}

\newtheorem{theorem}{Theorem}[section]
\newtheorem{corollary}[theorem]{Corollary}
\newtheorem{lemma}[theorem]{Lemma}

\newtheorem{conjecture}[theorem]{Conjecture}

\newtheorem{mainthm}[theorem]{Main Theorem}

\theoremstyle{definition}
\newtheorem{definition}[theorem]{Definition}
\newtheorem{remark}[theorem]{Remark}

\numberwithin{equation}{section}

\begin{document}


\baselineskip=17pt


\title[Perfect \textsc{Eisenstein} integers]{Perfect \textsc{Eisenstein} integers}

\author[J. C. Stumpenhusen]{Johann Christian Stumpenhusen}
\address{Universit\"at zu K\"oln\\
Mathematisches Institut\\
Weyertal 86 -- 90\\
50931 K\"oln}
\email{jstumpen@math.uni-koeln.de}

\date{June 2022}

\begin{abstract}
We generalise the sum-of-divisors-function $\sigma$ and evenness to the rings of integers of certain algebraic number fields. In particular, we present necessary and sufficient conditions for even \textsc{Eisenstein} integers to be (norm-)perfect based on the work of \textsc{McDaniel} \cite{Daniel} on \textsc{Gaussian} integers. Furthermore, some results concerning odd norm-perfect \textsc{Eisenstein} integers and the rings of integers of other cyclotomic fields are proven.
\end{abstract}



\maketitle

\section{Introduction}
The sum-of-divisors-function $\sigma$ is given for any natural number $n$ by
\begin{equation}\label{1stfdefsig}
    \sigma(n) = \sum_{d \, | \, n} d
\end{equation}
where the $d$'s are a natural numbers. It adds up all the positive divisors of $n$. Natural numbers satisfying
\begin{equation}\label{1stdefperf}
    \sigma(n) = 2n
\end{equation}
are called \textit{perfect numbers} and have already been studied by \textsc{Euclid}  (4\textsuperscript{th}/3\textsuperscript{rd} century BC), the ancient Greeks, and numerous mathematicians throughout history.

\subsection{Generalised $\sigma$-function and previous literature}
Back in the days of the ancient Greeks, it was already known to them that there is a significant difference between the even and odd perfect natural numbers. That is, they were solely able to find even numbers of that kind.

\begin{theorem}[\textsc{Euclid}--\textsc{Euler}]\label{EuclidEuler}
An even natural number is perfect if and only if it is of the form $2^{k-1}(2^k - 1)$ for some $k \in \mathbb{N}$ such that $2^k - 1$ is prime.
\end{theorem}

The \textit{if}-direction was proven by \textsc{Euclid} while the converse remained unknown until the Swiss mathematician \textsc{Leonard Euler} (1707 -- 1783 AD) proved it in 1747. The question of the even perfect rational integers may be regarded as completely solved, meanwhile the odd perfect rational integers are an open problem. Over the last few decades, there has been a significant progress towards the form which such a number can have. Noteworthy contributions have been made by \textsc{Nielsen} \cite{Nielsen} and \textsc{Ochem} and \textsc{Rao} \cite{OchRa} who provided an upper (based on the number of prime factors) and an unconditional lower bound for any odd perfect rational integer. Some other authors proved a lower bound for the amount of prime divisors of odd perfect rational integers. However, not a single such number has been found yet and it is uncertain whether there are any. When working with odd perfect integers, we will focus on the theorem \textsc{Euler} himself stated.

\begin{theorem}[\textsc{Euler}]\label{Euler}
Let $n$ be an odd perfect natural number. Then
$$n = p^kq^2$$
for some prime $p \in \mathbb{P}$ and odd natural number $q$. Moreover, $p \equiv k \equiv 1 \mod 4\mathbb{Z}$ and $\langle p, q \rangle = \mathbb{Z}$.
\end{theorem}

As the concept of divisibility is not restricted to $\mathbb{Z}$, the question arises how to generalise the $\sigma$-function to a class of rings which have similar properties: the rings of integers $\mathcal{O}_K$ for any number field $K$. It is well known that Equation \ref{1stfdefsig} may be rewritten as
\begin{equation}\label{2nddefsig}
    \sigma(n) = \prod_{k = 1}^l \frac{p_k^{e_k + 1} - 1}{p_k - 1}
\end{equation}
for a factorisation $n = \prod_{k = 1}^l p_k^{e_k}$ into natural prime numbers $p_k \in \mathbb{P}$. If $K$ has class number 1, we may choose an appropriate system of representatives among the primes elements of $\mathcal{O}_K$ (the \textit{positive} primes, denoted by $\mathbb{P}_K^+$) to define the \textit{generalised sum-of-divisors-function} for an integer $\alpha \in \mathcal{O}_K$ to be
\begin{equation}\label{sigK}
    \sigma_K(\alpha) = \prod_{k = 1}^l \frac{\pi_k^{e_k + 1} - 1}{\pi_k - 1}
\end{equation}
where $\alpha = \prod_{k = 1}^l \pi_k^{e_k}$ is a factorisation into prime elements $\pi_k \in \mathbb{P}_K^+$ of $\mathcal{O}_K$.

In the 20\textsuperscript{th} century, the concept of an even number was generalised to the ring of integers of some cyclotomic fields. A cyclotomic field is the splitting field of the polynomial $x^n - 1$ and generated by a primitive $n$-th root of unity $\zeta_n$. In this note, we will focus on the primitive roots of unity $\zeta_4 =: i$ and $\zeta_3 =: \omega$ especially.

\begin{definition}\label{defevenperfnormperf}
Let $K = \mathbb{Q}(\zeta_p)$ for some $p \in \mathbb{P}$ or $p = 4$. An element $\alpha \in \mathcal{O}_K$ is called \textit{even} if it is divisible by $1 - \zeta_p$. Otherwise, we call it \textit{odd}. An integer $\alpha \in \mathcal{O}_K$ is called
\begin{enumerate}
    \item \textit{perfect} if $\sigma_K(\alpha) = (1 - \zeta_p)\alpha$, and
    \item \textit{norm-perfect} if $N_K\left(\sigma_K(\alpha)\right) = N_K(1 - \zeta_p)N_K(\alpha)$.
\end{enumerate}
\end{definition}

\textsc{Spira} \cite{Spira} and \textsc{McDaniel} \cite{Daniel} did this generalisation explicitly for the \textsc{Gaussian} integers $\mathbb{Z}[i]$ and worked towards transferring the \textsc{Euclid}--\textsc{Euler} Theorem to this ring. The latter succeeded in doing so whereas \textsc{Parker}, \textsc{Rushall}, and \textsc{Hunt} \cite{Parker}, who tried the same for the \textsc{Eisenstein} integers $\mathbb{Z}[\omega]$, only managed to find an even norm-perfect integer.

A theorem regarding the odd norm-perfect \textsc{Gaussian} integers was published by \textsc{Ward} \cite{Ward}, using a similar approach as \textsc{Euler}. \textsc{Parker}, \textsc{Rushall}, and \textsc{Hunt} conjectured a form of odd norm-perfect \textsc{Eisenstein} integers, too.

\subsection{Results}
The first result refines the theorem by \textsc{Parker}, \textsc{Rushall}, and \textsc{Hunt}, presenting not only even norm-perfect but even perfect \textsc{Eisenstein} integers by using elegant inequalities by \textsc{Spira} \cite{Spira} and \textsc{McDaniel} \cite{Daniel}.

\begin{mainthm}\label{EuclidEulerEisenstein}
$\alpha \in \mathbb{Z}[\omega]$ is an even primitive norm-perfect number if and only if either
\begin{enumerate}
    \item $\alpha = \varepsilon(1 - \omega^2)^{k-1}[(1 - \omega^2)^k - 1]$ where $(1 - \omega^2)^k - 1$ is prime with a natural prime $k \equiv 1 \mod 12\mathbb{Z}$ or \label{EuclidEulerEisenstein1}
    \item $\alpha = \varepsilon(1 - \omega^2)^{k-1}\overline{[(1 - \omega^2)^k - 1]}$ where $(1 - \omega^2)^k - 1$ is prime with a natural prime $k \equiv -1 \mod 12\mathbb{Z}$ \label{EuclidEulerEisenstein2}
\end{enumerate}
and $\varepsilon$ is a unit in $\mathbb{Z}[\omega]$. Moreover, the even primitive perfect integers are exactly those $\alpha$ of item \ref{EuclidEulerEisenstein1} with $\varepsilon = - \omega$.
\end{mainthm}

We will also prove a theorem about the odd norm-perfect \textsc{Eisenstein} integers which will be followed by a generalisation to cyclotomic fields of higher degree.

\begin{theorem}\label{oddnormperfEisenstein}
Let $\alpha$ be an odd norm-perfect \textsc{Eisenstein} integer and define $P_j = \{\psi \in \mathbb{P}_{\mathbb{Q}(\omega)}^+: \psi \equiv j \mod 1 - \omega^2 \land \psi \, | \, \alpha\}$ for $j \in \{1,2\}$. $\alpha$ has to be of the form
$$\alpha = \varepsilon \psi_0^k\prod_{\psi_1 \in P_1, \, \psi_1 \neq \psi_0} \psi_1^{e_{\psi_1}}\prod_{\psi_2 \in P_2, \, \psi_2 \neq \psi_0} \psi_2^{e_{\psi_2}}$$
with either
\begin{enumerate}
    \item $\psi_0 \in P_1$, $k \equiv 2 \mod 3\mathbb{Z}$ or
    \item $\psi_0 \in P_2$, $k \equiv 1 \mod 2\mathbb{Z}$
\end{enumerate}
as well as $e_{\psi_1} \not \equiv 2 \mod 3\mathbb{Z}$ and $e_{\psi_2} \equiv 0 \mod 2\mathbb{Z}$ and $\varepsilon$ being a unit in either case.
\end{theorem}

\subsection{Notation}
Throughout this work, $K$ denotes an algebraic number field over the rational numbers, that is a finite algebraic extension of $\mathbb{Q}$. For any such $K$, the ring $\mathcal{O}_K$ contains all the elements of $K$ which satisfy a monic polynomial in $\mathbb{Z}[x]$, i.e. every $\alpha \in K$ such that $f(\alpha) = 0$ for some monic $f \in \mathbb{Z}[x]$. We call $\mathcal{O}_K$ the \textit{ring of integers} of $K$ and its elements \textit{algebraic integers} or \textit{integers} of $K$. In order to avoid confusion, the elements of $\mathbb{Z}$ are called \textit{rational} integers. For an element $\alpha \in K$, $N_K(\alpha)$ will denote its \textit{norm} in the extension $K/\mathbb{Q}$. The same holds true for an ideal $\mathfrak{a}$. As usual, we write $Cl(K)$ for the \textit{class group} of $K$.

Two elements $a, b$ in a ring $R$ are called \textit{associates} of each other if there is a unit $\varepsilon \in R$ such that $a = \varepsilon b$. We denote this relation by $a \simeq b$. If the ring is clarified by the context, we will write $\langle \alpha \rangle$ instead of $\alpha R$ for the ideal of $R$ generated by $\alpha$.

\section{The positive prime of minimal norm}\label{chapcycloperf}

We need to think about what the generalised $\sigma$-function is supposed to express. In the rational integers, it symbolises the accumulated size of all divisors of an integer. In order to conserve a somewhat constant ratio between the sizes of elements which we will deem interchangeable in some way, i.e. associates, we will restrict ourselves to the cyclotomic fields $\mathbb{Q}(\zeta_n)$ where $\zeta_n$ is an $n$-th root of unity. Their ring of integers is given by $\mathbb{Z}[\zeta_n]$.

\subsection{The positive primes}
The second property we need to address regarding the $\sigma$-function is that in the naive definition on the natural numbers (and their extension to the rational integers) we constrain the sum to positive divisors. However, the cyclotomic fields over $\mathbb{Q}$ are complex extensions and hence, we do not really have a sense of positivity, not even in $\mathbb{Z}[\zeta_n]$. Aiming towards creating something of that kind, we restrict the set of fields we want to work on even further: We only consider cyclotomic fields $K = \mathbb{Q}(\zeta_p)$ where $p \in \mathbb{P}$ is a natural prime or $p = 4$ and $\mathcal{O}_K$ is a UFD. We will call this set $\mathcal{R}$. The set of prime elements of an algebraic number field $K$ will be denoted by $\mathbb{P}_K$.

\begin{definition}
\label{defposcyclo} Let $K = \mathbb{Q}(\zeta_p) \in \mathcal{R}$. A prime $\psi \in \mathcal{O}_K$ is \textit{positive} if its angle with respect to the positive real line is smaller than any of its associates' and, if $p \neq 3, 4$, its absolute value in $\mathbb{C}$ is the smallest such that $|\psi| \geq \sqrt[p-1]{N_K(\psi)}$. An integer $\alpha$ of $K$ is \textit{positive} if it is the product of positive primes.  We denote the sets of positive primes or integers by $\mathbb{P}_K^+$ or $\mathcal{O}_K^+$, respectively.
\end{definition}

Of course, this settles which primitive $p$-th root of unity we consider in Definition \ref{defevenperfnormperf} if we assume $1 - \zeta_p$ to be the positive associate.

We note that the previous definition is coherent with $\mathbb{Q}$ being the field containing the second root of unity. The additional condition for $p \neq 3, 4$ is due to $\mathbb{Z}[\zeta_p]$ then containing an infinite group of units. It is also clearly different from the definition of being totally positive. Moreover, every non-zero integer has exactly one positive associate so that we can easily define a function that yields the same value for different associates just based on the positive integers of $K$. We point out that the choice of our positive set $\mathbb{P}_K^+$ is completely arbitrary, any representative system of the associate equivalence classes may be chosen. Since $\mathcal{O}_K$ is a UFD, the definition of $\sigma_K$ as in \ref{sigK} makes sense and transfers the $\sigma$-function to the integers of $K$.

The $\sigma$-function on the naturals numbers has the property that $\sigma(n) \geq n$ for all $n \in \mathbb{N}$. We will now investigate why the choice we made for $\mathcal{O}_K^+$ is advantageous. The following lemma was proven by \textsc{Spira} \cite{Spira}.

\begin{lemma}[\textsc{Spira}]
\label{SpiraIneq} Let $K \in \{\mathbb{Q}, \mathbb{Q}(\zeta_3), \mathbb{Q}(\zeta_4)\}$ and $\alpha \in \mathcal{O}_K$. We have
$$N_K\left(\frac{\alpha^{n+1} - 1}{\alpha - 1}\right) \geq N_K(\alpha^n)$$
if $\Re(\alpha) \geq 1$ and $\alpha \neq 0, 1$.
\end{lemma}

\subsection{Norms of primes}

Computing the discriminant of $K = \mathbb{Q}(\zeta_p)$ for an odd natural prime $p$ yields
$$\Delta_K = (-1)^{\frac{p - 1}{2}}p^{p-2}$$
and for $p = 4$, we have
$$\Delta_K = -4.$$
This implies that $p$ or $2$ (in the case of $p = 4$) is the only prime ramifying in $K$. Using the \textsc{Dedekind}--\textsc{Kummer} Theorem \cite[p. 46]{Orr}, we get
$$p\mathcal{O}_K = \langle 1 - \zeta_p\rangle^{p-1}$$
if $p$ is odd and
$$2\mathbb{Z}[i] = \langle 1 + i\rangle^2$$
if $p = 4$ where $\langle 1 - \zeta_p\rangle$ is, indeed, a prime ideal of $\mathcal{O}_K$, so $p$ (or $2$ resp.) totally ramifies. We also compute $N_K(1 - \zeta_p) = p$ if $p$ is odd and $N_K(1 + i) = 2$ for $p = 4$.

For any other prime, we may use the following lemma for computing the norm of the prime ideals above it after recalling that the norm $N(\mathfrak{a})$ of an ideal $\mathfrak{a}$ is equal to the norm of its generator if $\mathcal{O}_K$ is a PID.

\begin{lemma}
\label{lemminimalprime} Let $p, q \in \mathbb{P}$, $p \neq q$, be positive rational primes and $\mathfrak{Q} \supset \langle q \rangle$ a prime ideal in $K = \mathbb{Q}(\zeta_p)$. Then $N_K(\mathfrak{Q}) = q^f$ where $f$ is the order of $q$ in $\left(\nicefrac{\mathbb{Z}}{p\mathbb{Z}}\right)^*$.
\end{lemma}

\textit{Proof.} Because of \textsc{Dekekind}--\textsc{Kummer}, we consider the equation
$$x^p - 1 \equiv 0 \mod q\mathbb{Z}$$
since any  of its irreducible factors other than $x - 1$ is an irreducible factor of $g := \sum_{k = 0}^{p-1} x^k \mod q\mathbb{Z}$ which is the reduction of the minimal polynomial of $\zeta_p$. Let $F$ be the smallest extension of $\mathbb{F}_q$ such that it contains all the roots of $x^p - 1$, say $F = \mathbb{F}_{q^f}$. By \textsc{Lagrange}'s Theorem about the order of an element, it is $p \, | \, q^f - 1$ and hence, $q^f \equiv 1 \mod p\mathbb{Z}$. As $F$ is the smallest such extension, we also deduce that $f$ is the smallest positive integer satisfying that congruence condition. Thus, each irreducible factor of $g$ has degree $f$ and it follows $N_K(\mathfrak{Q}) = q^f$. \hfill $\Box$\\

For any positive rational prime $q < p$, we have $q > 1$ and thus $q^f > p$ if $q^f \equiv 1 \mod p$. Hence, we may say that $1 - \zeta_p$ is the non-unit of minimal norm because the norm is a multiplicative function. Moreover, we see that the norm only takes values being $\equiv 0, 1 \mod p\mathbb{Z}$ because each prime element has norm satisfying that congruence condition as we saw above.

\subsection{Generalised \textsc{Mersenne} numbers}

Lastly, we want to generalise the \textsc{Mersenne} numbers.

\begin{definition}\label{defMersenne}
The natural numbers $2^k - 1$ for $k \in \mathbb{N}$ are called \textsc{Mersenne} numbers.
\end{definition}

Due to their exponential distribution, \textsc{Mersenne} primes are still among the largest known prime numbers. A necessary condition for such a number to be prime is that $k$ is prime itself. However, this is not sufficient as we may see by the example
$$2^{11} - 1 = 2047 = 23 \cdot 89.$$
The Theorem \ref{EuclidEuler} solved the question of the form of even perfect natural numbers and proved that there is a bijection between these numbers and \textsc{Mersenne} primes. It does not make any statement about the amount of such numbers. Presently, the \textit{Great Internet Mersenne Prime Search} \cite{GIMPS}, abbreviated GIMPS, is working on testing several large \textsc{Mersenne} numbers on primality.

\begin{definition}\label{defmersennecyclo}
Let $K = \mathbb{Q}(\zeta_p) \in \mathcal{R}$. We call the numbers $(1 - \zeta_p)^k - 1$ with $k \in \mathbb{N}$ the \textsc{Mersenne} numbers of $K$.
\end{definition}

It is important to note that, contrarily to the case of the rational integers, there is no bijection from the perfect \textsc{Gaussian} or \textsc{Eisenstein} integers to the \textsc{Mersenne} primes of the respective field but only to a subset thereof. However, there is a result concerning the rational \textsc{Mersenne} numbers which may be generalised to any field $K \in \mathcal{R}$.

\begin{lemma}\label{expomersennecyclo}
If $(1 - \zeta_p)^k - 1$ with $k \in \mathbb{N}_{\geq 2}$ is prime, then $k$ is a rational prime.
\end{lemma}

\textit{Proof.} Suppose $k$ is not prime, i.e. $k = mn$ with $m, n \in \mathbb{N}_{\geq 2}$. Then
\begin{align*}
    (1 - \zeta_p)^k - 1 &= [(1 - \zeta_p)^m]^n - 1\\
    &= [(1 - \zeta_p)^m - 1]\sum_{j = 0}^{n-1} (1 - \zeta_p)^{mj}.
\end{align*}
As $m \geq 2$, we have that the left factor is not a unit. Moreover, as $n \geq 2$, we have $m < k$, so $(1 - \zeta_p)^m - 1 \not \simeq (1 - \zeta_p)^k - 1$.\hfill $\Box$\\

In the case of the \textsc{Gaussian} and \textsc{Eisenstein} \textsc{Mersenne} numbers, \textsc{Berrizbeitia} and \textsc{Iskra} \cite{BerrizbeitiaIskra} presented primality tests similar to the \textsc{Lucas}--\textsc{Lehmer} test for the rational \textsc{Mersenne} numbers.

\section{Perfect numbers in the cyclotomic fields}

For the sake of completeness, we transfer one more definition to the ring of integer of a field $K \in \mathcal{R}$.

\begin{definition}
Let $K \in \mathcal{R}$. An integer $\alpha \in \mathcal{O}_K$ is called
\begin{enumerate}
    \item \textit{abundant} if $N_K\left(\sigma_K(\alpha)\right) > N_K(1 - \zeta_p)N_K(\alpha)$ or
    \item \textit{deficient} if $N_K\left(\sigma_K(\alpha)\right) < N_K(1 - \zeta_p)N_K(\alpha)$.
\end{enumerate}
\end{definition}

\subsection{\textsc{Gaussian} integers}

One of the most commonly known algebraic extensions of the rationals $\mathbb{Q}$ is $\mathbb{Q}(i)$. The norm function is given by $N_{\mathbb{Q}(i)}(a + bi) = a^2 + b^2$ and its set of positive primes is $\mathbb{P}_{\mathbb{Q}(i)}^+ = \{a + bi \in \mathbb{Z}[i]: a > 0, b \geq 0\} \cap \mathbb{P}_{\mathbb{Q}(i)}$. Thus, the minimal prime is $1 + i$.

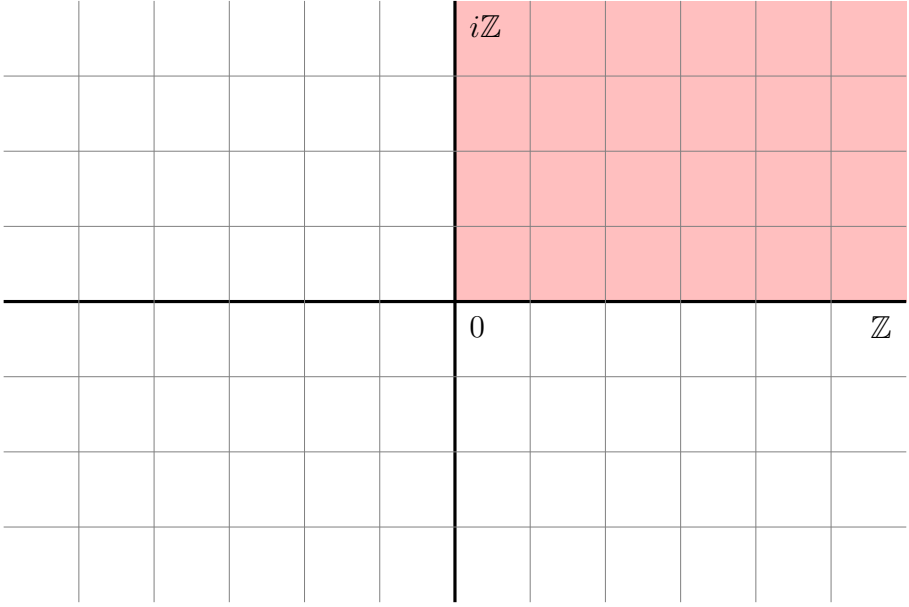
\begin{figure}[H]
    \centering
    \begin{tikzpicture}
    
    \fill[red!25] (0,0) -- (6,0) -- (6,4) -- (0,4) -- cycle;
    
    \node[below right=2pt of {(0,0)},fill=white]{$0$};
    \node[below right=2pt of {(0,4)}]{$i\mathbb{Z}$};
    \node[below left=2pt of {(6,0)}]{$\mathbb{Z}$};
    
    \draw[black, very thick] (-6, 0) -- (6, 0);
    \draw[black, very thick] (0, -4) -- (0, 4);
    
    \draw[gray, ultra thin] (-6, 1) -- (6, 1);
    \draw[gray, ultra thin] (-6, 2) -- (6, 2);
    \draw[gray, ultra thin] (-6, 3) -- (6, 3);
    \draw[gray, ultra thin] (-6, -1) -- (6, -1);
    \draw[gray, ultra thin] (-6, -2) -- (6, -2);
    \draw[gray, ultra thin] (-6, -3) -- (6, -3);
    
    \draw[gray, ultra thin] (5, 4) -- (5, -4);
    \draw[gray, ultra thin] (4, 4) -- (4, -4);
    \draw[gray, ultra thin] (3, 4) -- (3, -4);
    \draw[gray, ultra thin] (2, 4) -- (2, -4);
    \draw[gray, ultra thin] (1, 4) -- (1, -4);
    \draw[gray, ultra thin] (-1, 4) -- (-1, -4);
    \draw[gray, ultra thin] (-2, 4) -- (-2, -4);
    \draw[gray, ultra thin] (-3, 4) -- (-3, -4);
    \draw[gray, ultra thin] (-4, 4) -- (-4, -4);
    \draw[gray, ultra thin] (-5, 4) -- (-5, -4);
    \end{tikzpicture}
    \label{Gausspic}
    \caption{The \textsc{Gaussian} integers near the origin, each represented by an intersection of two lines. The positive primes lie in the shaded domain, excluding the left boundary.}
\end{figure}

If we pay a close look to the behaviour of odd and even numbers in $\mathbb{Z}[i]$, we see that, due to $N_{\mathbb{Q}(i)}(1 + i) = 2$, it is
$$\nicefrac{\mathbb{Z}[i]}{\langle 1 + i \rangle} \cong \nicefrac{\mathbb{Z}}{2\mathbb{Z}}.$$
Hence, there is only one congruence class of odd \textsc{Gaussian} integers and the computations modulo $(1 + i)\mathbb{Z}[i]$ are similar to those modulo $2\mathbb{Z}$ in the rational integers. We may therefore distinguish between odd and even perfect \textsc{Gaussian} integers.

This subsection presents the main ideas used by \textsc{McDaniel} \cite{Daniel} in order to transfer the \textsc{Euclid}--\textsc{Euler} Theorem to the \textsc{Gaussian} integers. A crucial step is the improvement of \textsc{Spira}'s inequality in Lemma \ref{SpiraIneq}.

\begin{lemma}[\textsc{McDaniel} {{\cite[p. 138]{Daniel}}}] \label{DanielIneq}
Let $\alpha \in \mathbb{Z}[i]$ and $|\alpha| \geq \sqrt{5}$. Then
$$N_{\mathbb{Q}(i)}\left(\frac{\alpha^{n+1} - 1}{\alpha - 1}\right) > N_{\mathbb{Q}(i)}(\alpha^n)\left(1 + \frac{2\Re(\alpha) - \frac{7}{5}}{N_{\mathbb{Q}(i)}(\alpha)}\right)$$
for all $n \in \mathbb{N}_{>0}$.
\end{lemma}

The proof consists mostly of computations. However, using the definition of $\sigma_{\mathbb{Q}(i)}$ and the multiplicativity of the norm function, this lemma yields a very useful corollary.

\begin{corollary}\label{DanielIneqCoro}
Let $\psi$ be an odd positive \textsc{Gaussian} prime. Then
$$N_{\mathbb{Q}(i)}\left(\frac{\sigma_{\mathbb{Q}(i)}(\psi^n)}{\psi^n}\right) > 1 + \frac{2\Re(\psi) - \frac{7}{5}}{N_{\mathbb{Q}(i)}(\psi)}$$
for all $n \in \mathbb{N}_{>0}$.
\end{corollary}

In order to use the previous lemma to prove this corollary, we point out that any odd positive \textsc{Gaussian} prime has an absolute value of at least $\sqrt{5}$.

When taking a closer look at $N_{\mathbb{Q}(i)}(\alpha)$ for an arbitrary $\alpha \in \mathbb{Z}[i]$, we get an inequality which looks quite similar to one we are familiar with from the rational integers.

\begin{lemma}\label{sigcycloIneq}
Let $\alpha \in \mathbb{Z}[i]$. Then $N_{\mathbb{Q}(i)}\left(\sigma_{\mathbb{Q}(i)}(\alpha)\right) \geq N_{\mathbb{Q}(i)}(\alpha)$.
\end{lemma}

\textit{Proof.} The norm on $\mathbb{Z}[i]$ is multiplicative and any positive \textsc{Gaussian} prime satisfies Lemma \ref{SpiraIneq}.\hfill $\Box$\\

In the rational integers, the \textsc{Eulcid}--\textsc{Euler} Theorem is restricted to even perfect numbers and, similarly, \textsc{McDaniel}'s theorem only answers the question about even (norm-) perfect \textsc{Gaussian} integers. Corollary \ref{DanielIneqCoro} and Lemma \ref{sigcycloIneq} together yield some lower bounds on the norm of an odd prime divisor of an even norm-perfect \textsc{Gaussian} integer. These are joint by \textsc{McDaniel} with a theorem by \textsc{Spira}.

\begin{lemma}[\textsc{Spira} {{\cite[p. 123]{Spira}}}] Let $\eta = (1 + i)^{k-1}\mu$ be a norm-perfect \textsc{Gaussian} integer with $k \in \mathbb{N}_{\geq 2}$ and $\mu$ odd. Then $k \equiv 0, \pm 1 \mod 12\mathbb{Z}$.
\end{lemma}

Before we move on, there is still one other definition which was used in the main theorem of this note.

\begin{definition}\label{defprimnormperf}
Let $\eta$ be a (norm-)perfect \textsc{Gaussian} integer. $\eta$ is called \textit{primitive} if there is no $\theta \in \mathbb{Z}[i]$ such that $\theta$ is (norm-)perfect, $\theta \, | \, \eta$ and $\theta \not \simeq \eta$.
\end{definition}

This deals with the fact that Lemma \ref{sigcycloIneq} does not provide a strict inequality.

\begin{theorem}[\textsc{McDaniel} {{\cite[p. 137]{Daniel}}}]
\label{Danielmain} Let $M_p = (1 + i)^p - 1$ be a \textsc{Gaussian Mersenne} prime and $\varepsilon$ a unit in $\mathbb{Z}[i]$.
\begin{enumerate}
    \item If $p \equiv 1 \mod 8\mathbb{Z}$, then $\eta = \varepsilon(1 + i)^{p-1}M_p$ is a primitive norm-perfect \textsc{Gaussian} integer.
    \item If $p \equiv -1 \mod 8\mathbb{Z}$, then $\eta = \varepsilon(1 + i)^{p-1}\overline{M_p}$ is a primitive norm-perfect \textsc{Gaussian} integer.
\end{enumerate}
Conversely, if $\eta$ is an even norm-perfect \textsc{Gaussian} integer, then, for some unit $\varepsilon$, there is either
\begin{enumerate}
    \item a rational prime $p \equiv 1 \mod 8\mathbb{Z}$ such that $\eta = \varepsilon(1 + i)^{p-1}M_p$ or
    \item a rational prime $p \equiv -1 \mod 8\mathbb{Z}$ such that $\eta = \varepsilon(1 + i)^{p-1}\overline{M_p}$
\end{enumerate}
where $M_p$ is a \textsc{Gaussian Mersenne} prime. Moreover, \textit{norm-perfect} may be substituted by \textit{perfect} if we only consider the first bullet in each part and replace $\varepsilon$ by $-i$.
\end{theorem}

This settles the even (norm-)perfect \textsc{Gaussian} integers. So far, there has not been much work about their odd counterparts but \textsc{Ward} proved a theorem on the form of such an integer.

\begin{theorem}[\textsc{Ward} {{\cite[p. 2]{Ward}}}] \label{Ward}
Let $\eta$ be an odd norm-perfect \textsc{Gaussian} integer. Then it is of the form
$$\eta = \psi^k\rho^2$$
for an odd $k \in \mathbb{N}$, an odd \textsc{Gaussian} prime $\psi$, and an odd \textsc{Gaussian} integer $\rho$.
\end{theorem}

Surprisingly and in contrast to the rational integers where all primes are deficient, \textsc{Ward} is able to present an odd positive norm-perfect \textsc{Gaussian} integer: $2 + i$ which is also prime. We quickly check that
$$N_{\mathbb{Q}(i)}(\sigma_{\mathbb{Q}(i)}(2 + i)) = N_{\mathbb{Q}(i)}(3 + i) = 3^2 + 1^2 = 10 = 2 \cdot 5 = 2 N_{\mathbb{Q}(i)}(2+i).$$
However, this is the only norm-perfect \textsc{Gaussian} prime which can be seen by evaluating the equation
$$N_{\mathbb{Q}(i)}(\psi + 1) = N_{\mathbb{Q}(i)}\left(\sigma_{\mathbb{Q}(i)}(\psi)\right) = 2N_{\mathbb{Q}(i)}(\psi)$$
as $\sigma_{\mathbb{Q}(i)}(\psi) = \psi + 1$ for any positive prime $\psi$. If we choose our set of positive primes $\mathbb{P}_{\mathbb{Q}(i)}^+$ differently, such that $2 - i$ is positive instead of $1 + 2i$, it will be another odd positive norm-perfect prime. This underlines the fact that our set of (norm-)perfect integers depends on the set $\mathbb{P}_K^+$ for $K \in \mathcal{R}$.

\subsection{\textsc{Eisenstein} integers}

The degree of the field extension $[\mathbb{Q}(\zeta_n) : \mathbb{Q}]$ is $\varphi(n)$ where $\varphi$ is \textsc{Euler}'s totient function. There are exactly two solutions to $\varphi(n) = 2$ which are $n \in \{3,4\}$. We dealt with $n = 4$ in the previous subsection but less has been done in the case of $n = 3$ so far.

\begin{definition}\label{defomega}
We define
$$\omega := \exp\left(\frac{2\pi i}{3}\right),$$
the primitive third root of unity in the upper complex half-plane.
\end{definition}

The extensions $\mathbb{Q}(i)/\mathbb{Q}$ and $\mathbb{Q}(\omega)/\mathbb{Q}$ have a lot in common. We recapitulate some of the properties of the latter, being a cyclotomic field of degree 2 over the rationals. The norm function is given by $N_{\mathbb{Q}(\omega)}(a + b\omega) = a^2 - ab + b^2$ und its set of positive primes is $\mathbb{P}_{\mathbb{Q}(\omega)}^+ = \{a + b\omega \in \mathbb{Z}[\omega]: a > b \geq 0\} \cap \mathbb{P}_{\mathbb{Q}(\omega)}$. Thus, the minimal prime is $2 + \omega = 1 - \omega^2$.

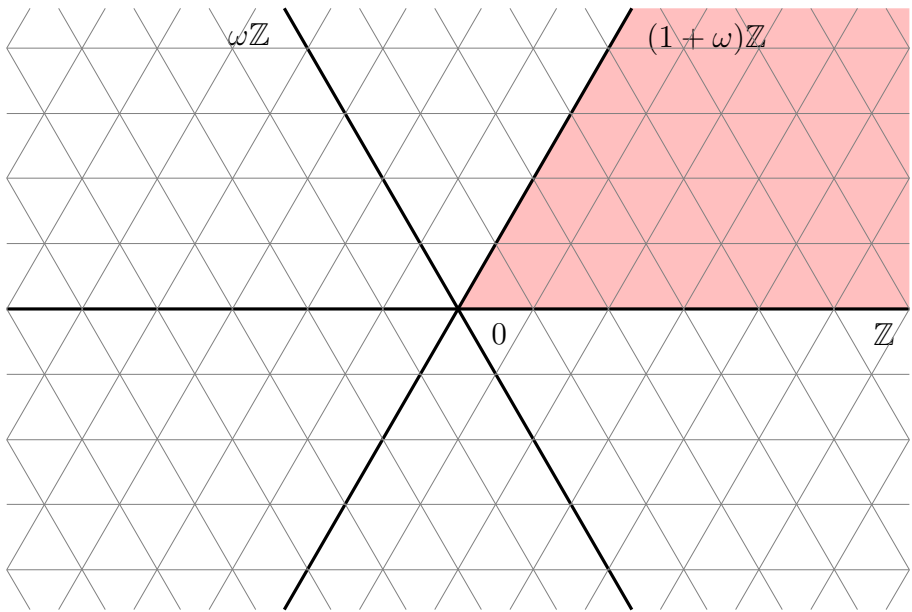
\begin{figure}[H]
    \centering
    \begin{tikzpicture}
    \fill[red!25] (0,0) -- (6,0) -- (6,4) -- (2.31,4) -- cycle;
    
    \node[below right=2pt of {(0.25,0)},fill=white]{$0$};
    \node[below left=2pt of {(6,0)},fill=white]{$\mathbb{Z}$};
    \node[below right=2pt of {(2.31,4)},fill=red!25]{$(1 + \omega)\mathbb{Z}$};
    \node[below left=2pt of {(-2.31,4)},fill=white]{$\omega\mathbb{Z}$};
    
    \draw[black, very thick] (-6, 0) -- (6, 0);
    \draw[black, very thick] (-2.31, -4) -- (2.31, 4);
    \draw[black, very thick] (2.31, -4) -- (-2.31, 4);
    
    \draw[gray, ultra thin] (3.31, -4) -- (-1.31, 4);
    \draw[gray, ultra thin] (4.31, -4) -- (-0.31, 4);
    \draw[gray, ultra thin] (5.31, -4) -- (0.69, 4);
    \draw[gray, ultra thin] (6, -3.47) -- (1.69, 4);
    \draw[gray, ultra thin] (6, -1.74) -- (2.69, 4);
    \draw[gray, ultra thin] (6, 0) -- (3.69, 4);
    \draw[gray, ultra thin] (6, 1.74) -- (4.69, 4);
    \draw[gray, ultra thin] (6, 3.47) -- (5.69, 4);
    
    \draw[gray, ultra thin] (-3.31, -4) -- (1.31, 4);
    \draw[gray, ultra thin] (-4.31, -4) -- (0.31, 4);
    \draw[gray, ultra thin] (-5.31, -4) -- (-0.69, 4);
    \draw[gray, ultra thin] (-6, -3.47) -- (-1.69, 4);
    \draw[gray, ultra thin] (-6, -1.74) -- (-2.69, 4);
    \draw[gray, ultra thin] (-6, 0) -- (-3.69, 4);
    \draw[gray, ultra thin] (-6, 1.74) -- (-4.69, 4);
    \draw[gray, ultra thin] (-6, 3.47) -- (-5.69, 4);
    
    \draw[gray, ultra thin] (-3.31, 4) -- (1.31, -4);
    \draw[gray, ultra thin] (-4.31, 4) -- (0.31, -4);
    \draw[gray, ultra thin] (-5.31, 4) -- (-0.69, -4);
    \draw[gray, ultra thin] (-6, 3.47) -- (-1.69, -4);
    \draw[gray, ultra thin] (-6, 1.74) -- (-2.69, -4);
    \draw[gray, ultra thin] (-6, 0) -- (-3.69, -4);
    \draw[gray, ultra thin] (-6, -1.74) -- (-4.69, -4);
    \draw[gray, ultra thin] (-6, -3.47) -- (-5.69, -4);
    
    \draw[gray, ultra thin] (3.31, 4) -- (-1.31, -4);
    \draw[gray, ultra thin] (4.31, 4) -- (-0.31, -4);
    \draw[gray, ultra thin] (5.31, 4) -- (0.69, -4);
    \draw[gray, ultra thin] (6, 3.47) -- (1.69, -4);
    \draw[gray, ultra thin] (6, 1.74) -- (2.69, -4);
    \draw[gray, ultra thin] (6, 0) -- (3.69, -4);
    \draw[gray, ultra thin] (6, -1.74) -- (4.69, -4);
    \draw[gray, ultra thin] (6, -3.47) -- (5.69, -4);
    
    \draw[gray, ultra thin] (6, 0.87) -- (-6, 0.87);
    \draw[gray, ultra thin] (6, 1.74) -- (-6, 1.74);
    \draw[gray, ultra thin] (6, 2.6) -- (-6, 2.6);
    \draw[gray, ultra thin] (6, 3.47) -- (-6, 3.47);
    \draw[gray, ultra thin] (6, -0.87) -- (-6, -0.87);
    \draw[gray, ultra thin] (6, -1.74) -- (-6, -1.74);
    \draw[gray, ultra thin] (6, -2.6) -- (-6, -2.6);
    \draw[gray, ultra thin] (6, -3.47) -- (-6, -3.47);
    \end{tikzpicture}
    \caption{The \textsc{Eisenstein} integers near the origin, each represented by an intersection of three lines. The positive primes lie in the shaded domain, excluding the upper boundary.}
    \label{Eisenpic}
\end{figure}

A main difference to $\mathbb{Z}[i]$ and $\mathbb{Z}$ directly impacting the computations we want to do in $\mathbb{Z}[\omega]$ is that, since $N_{\mathbb{Q}(\omega)}(1 - \omega^2) = 3$,
$$\nicefrac{\mathbb{Z}[\omega]}{\langle 1 - \omega^2\rangle} \cong \nicefrac{\mathbb{Z}}{3\mathbb{Z}}.$$
So there are two congruence classes of odd \textsc{Eisenstein} integers, namely those congruent to 1 or 2 modulo $(1 - \omega^2)\mathbb{Z}[\omega]$.

We will also copy the definition of a primitive (norm-)perfect number.

\begin{definition}\label{defprimperfeisen}
Let $\alpha$ be a (norm-)perfect \textsc{Eisenstein} integer. $\alpha$ is called \textit{primitive} if there is no $\beta \in \mathbb{Z}[\omega]$ such that $\beta$ is (norm-)perfect, $\beta \, | \, \alpha$ and $\beta \not \simeq \alpha$.
\end{definition}

In 2016, \textsc{Parker}, \textsc{Rushall}, and \textsc{Hunt} published a short paper \cite{Parker} transferring some of \textsc{McDaniel}'s results to the \textsc{Eisenstein} integers. Even though they use a similar definition of $\sigma_{\mathbb{Q}(\omega)}$, they do not introduce the concept of positive \textsc{Eisenstein} integers. So for simplicity, they choose $1 - \omega$ to be their minimal prime while the relevant---in our sense of definition \ref{defposcyclo} \textit{positive}---associate is still chosen to be $1 - \omega^2$. However, they present a sufficient condition for an even norm-perfect \textsc{Eisenstein} integer. The theorem is cited in its original form and therefore uses \textsc{Parker} et al.'s definitions.

\begin{theorem}[\textsc{Parker}, \textsc{Rushall}, \textsc{Hunt} {{\cite[p. 10]{Parker}}}] \label{Park1}
Given any rational integer $k > 1$, if $(1 - \omega)^k - 1$ is an \textsc{Eisenstein Mersenne} prime and if $k \equiv 11 \mod 12\mathbb{Z}$, then $\alpha = (1 - \omega)^{k-1}[(1 - \omega)^k - 1]$ is an even norm-perfect \textsc{Eisenstein} integer.
\end{theorem}

Oddly enough, \textsc{Parker} et al. do not refine their computations in this theorem's proof to show that---with their definitions---a particular associate of $\alpha$ is indeed perfect. We will adjust the theorem to our own definitions in the next subsection and close this gap. The main aim of this part of the note will be to also prove a converse statement, i.e. to prove the counterpart of the \textsc{Euclid}--\textsc{Euler} Theorem for \textsc{Eisenstein} integers. The inspiration for this is given by \textsc{McDaniel}'s Theorem \ref{Danielmain}.

\textsc{Parker} et al. also state a conjecture about the form of an odd norm-perfect \textsc{Eisenstein} integer.

\begin{conjecture}[{{\cite[p. 11]{Parker}}}]
\label{conPark} Any odd norm-perfect \textsc{Eisenstein} integer has to be of the form $\alpha = \psi^k\gamma^3$ where $\psi$ is an odd \textsc{Eisenstein} prime, $k \equiv 2 \mod 3\mathbb{Z}$ a rational integer and $\gamma$ an odd \textsc{Eisenstein} integer coprime to $\psi$.
\end{conjecture}

So far, there has not been any proof for this conjecture to be true but we will prove another form in the next-again subsection, presenting this form as a mere subcase.

\subsection{Even perfect Eisenstein integers}

The aim of this section is to prove the main theorem of this note, transferring the \textsc{Euclid}--\textsc{Euler} Theorem to the \textsc{Eisenstein} integers. Notice the similarity to the results by \textsc{McDaniel}. The significantly more difficult part is proving the \textit{only-if}-direction in either part of the theorem since computing the value of the $\sigma_{\mathbb{Q}(\omega)}$-function of a given element is fairly easy. To this end, we transfer two  of \textsc{Spira}'s lemmata.

\begin{lemma}\label{sigIneqeisen}
Let $\alpha \in \mathbb{Z}[\omega]$. Then $N_{\mathbb{Q}(\omega)}\left(\sigma_{\mathbb{Q}(\omega)}(\alpha)\right) \geq N_{\mathbb{Q}(\omega)}(\alpha)$.
\end{lemma}

\textit{Proof.} The norm is multiplicative and any positive \textsc{Eisenstein} prime satisfies Lemma \ref{SpiraIneq}.\hfill $\Box$\\

In particular, this lemma tells us that any multiple of an abundant number is also abundant. This helps us say something about the exponent $k$ of an even norm-perfect number.

\begin{lemma}
\label{exponormperfeisen}
Let $\alpha = (1 - \omega^2)^{k-1}\mu$ be an even norm-perfect \textsc{Eisenstein} integer with $\mu$ odd and a rational $k \geq 2$. Then $k \equiv \pm 2, \pm 1, 0 \mod 12\mathbb{Z}$.
\end{lemma}

\textit{Proof.} Recall that the norm and $\sigma_{\mathbb{Q}(\omega)}$-function are multiplicative and thus, so is their composition. We decompose
\begin{align*}
    3N_{\mathbb{Q}(\omega)}(\alpha) &= N_{\mathbb{Q}(\omega)}\left(\sigma_{\mathbb{Q}(\omega)}(\alpha)\right)\\
    &= N_{\mathbb{Q}(\omega)}\left(\sigma_{\mathbb{Q}(\omega)}((1 - \omega^2)^{k-1})\right)N_{\mathbb{Q}(\omega)}\left(\sigma_{\mathbb{Q}(\omega)}(\mu)\right).
\end{align*}
We are interested in the left factor and continue our computations.
\begin{align*}
    N_{\mathbb{Q}(\omega)}(\sigma_{\mathbb{Q}(\omega)}((1 - \omega^2)^{k-1})) &= N_{\mathbb{Q}(\omega)}\left(\frac{(1 - \omega^2)^k - 1}{1 - \omega^2 - 1}\right)\\
    &= N_{\mathbb{Q}(\omega)}\left(-\omega((1 - \omega^2)^k - 1)\right)\\
    &= 3^k - 2\Re\left((1 - \omega^2)^k\right) + 1.
\end{align*}
For $k \not \equiv \pm 2, \pm 1, 0 \mod 12\mathbb{Z}$, the middle term is less than or equal to 0, so the whole expression is greater than $3^k = 3N_{\mathbb{Q}(\omega)}\left((1 - \omega^2)^{k-1}\right)$. Using Lemma \ref{sigIneqeisen}, this would imply that $\alpha$ is abundant. Therefore, $k \equiv \pm 2, \pm 1, 0 \mod 12\mathbb{Z}$.\hfill $\Box$\\

This considerably reduces the possible residue classes of $k \mod 12\mathbb{Z}$. Now, we will extend Lemma \ref{DanielIneq} and its corollary \ref{DanielIneqCoro} to the \textsc{Eisenstein} integers.

\begin{lemma}
Let $\psi$ be an odd positive \textsc{Eisenstein} prime. Then
$$N_{\mathbb{Q}(\omega)}\left(\frac{\sigma_{\mathbb{Q}(\omega)}(\psi^n)}{\psi^n}\right) > 1 + \frac{2\Re(\psi) - \frac{7}{5}}{N_{\mathbb{Q}(\omega)}(\psi)}$$
for all $n \in \mathbb{N}_{>0}$.
\end{lemma}

\textit{Proof.} Firstly, we remark that the computations in the proof of Lemma \ref{DanielIneq} \cite{Daniel} are based on the norm on $\mathbb{Z}[i]$. Since the automorphism groups of both $\mathbb{Q}(i)$ and $\mathbb{Q}(\omega)$ only consist of the identity and the complex conjugation and therefore have a similar image, we may easily transfer that lemma to the \textsc{Eisenstein} integers for $|\psi| > \sqrt{5}$.

We are left to examine the inequality for all odd positive primes $\psi$ with absolute value less than $\sqrt{5}$. There is only one such prime, namely 2. However, we have
\begin{align*}
    N_{\mathbb{Q}(\omega)}\left(\frac{2^{n+1} - 1}{2 - 1}\right) &= N_{\mathbb{Q}(\omega)}(2^n) N_{\mathbb{Q}(\omega)}\left(2 - \frac{1}{2^n}\right)\\
    &= N_{\mathbb{Q}(\omega)}(2^n)\left(4 - \frac{1}{2^{n-2}} + \frac{1}{2^{2n}}\right)\\
    &> N_{\mathbb{Q}(\omega)}(2^n) \left(1 + \frac{\frac{13}{5}}{4}\right)\\
    &= N_{\mathbb{Q}(\omega)}(2^n) \left(1 + \frac{2\Re(2) - \frac{7}{5}}{N_{\mathbb{Q}(\omega)}(2)}\right)
\end{align*}
for all $n \in \mathbb{N}_{>0}$. We deduce the desired inequality by dividing both sides by $N_{\mathbb{Q}(\omega)}(2^n)$. \hfill$\Box$\\

We move on to find some estimate about the norm of an odd divisor of an even norm-perfect \textsc{Eisenstein} integer, so that we can minimise the number of possible prime divisors. This technique is inspired by \textsc{McDaniel}'s work but a bit more straight forward, since the setting of the \textsc{Eisenstein} integers lets us drop one extra lemma that was crucial in \textsc{McDaniel}'s proof.

\begin{lemma}\label{Ineqodddiv}
Let $\alpha = (1 - \omega^2)^{k-1}\mu$ be a norm-perfect \textsc{Eisenstein} integer with $\mu$ odd and a rational $k \geq 2$. Let $\psi$ be a positive prime divisor of $\mu$. Then
$$N(\psi) > \frac{13}{5} \cdot \frac{N_{\mathbb{Q}(\omega)}\left((1 - \omega^2)^k - 1\right)}{3^k - N_{\mathbb{Q}(\omega)}\left((1 - \omega^2)^k - 1\right)}.$$
\end{lemma}

\textit{Proof.} Let $e$ be the maximal exponent such that $\psi^e \, | \, \mu$. Using the fact that $\alpha$ is norm-perfect and Lemma \ref{sigIneqeisen}, we get
$$1 = \frac{N_{\mathbb{Q}(\omega)}\left(\sigma_{\mathbb{Q}(\omega)}(\alpha)\right)}{3N_{\mathbb{Q}(\omega)}(\alpha)} \geq \frac{N_{\mathbb{Q}(\omega)}\left((1 - \omega^2)^k - 1\right)N_{\mathbb{Q}(\omega)}\left(\sigma_{\mathbb{Q}(\omega)}(\psi^e)\right)}{3^kN_{\mathbb{Q}(\omega)}(\psi^e)}.$$
The previous lemma then tells us that
$$\frac{N_{\mathbb{Q}(\omega)}\left((1 - \omega^2)^k - 1\right)N_{\mathbb{Q}(\omega)}\left(\sigma_{\mathbb{Q}(\omega)}(\psi^e)\right)}{3^kN_{\mathbb{Q}(\omega)}(\psi^e)}$$
$$> \frac{N_{\mathbb{Q}(\omega)}\left(((1 - \omega^2)^k - 1\right)\left(N_{\mathbb{Q}(\omega)}(\psi) + 2\Re(\psi) - \frac{7}{5}\right)}{3^k\left(N_{\mathbb{Q}(\omega)}(\psi)\right)}.$$
Solving this for $N_{\mathbb{Q}(\omega)}(\psi)$ yields
$$N_{\mathbb{Q}(\omega)}(\psi) > (2\Re(\psi) - \frac{7}{5})\frac{N_{\mathbb{Q}(\omega)}\left((1 - \omega^2)^k - 1\right)}{3^k - N_{\mathbb{Q}(\omega)}\left((1 - \omega^2)^k - 1\right)}.$$
We have a quick look at the smallest odd positive primes of $\mathbb{Z}[\omega]$ so that we can estimate $2\Re(\psi)$. In the following table, $p$ is the rational prime such that $\psi$ lies above $p$.
\begin{table}[H]
    \centering
    $$\begin{array}{c|c|c}
        p & \psi &  \Re(\psi) \\
        \hline
        2 & 2 & 2 \\
        5 & 5 & 5 \\
        7 & 2 - \omega^2 & \frac{5}{2} \\
        7 & 1 - 2\omega^2 & 2 \\
        11 & 11 & 11 
    \end{array}$$
    \caption{Real parts of the odd positive primes of the least norm.}
    \label{RealPartsOddEisen}
\end{table}
Thus, $\Re(\psi) \geq 2$ and we are done since $2 \cdot 2 - \frac{7}{5} = \frac{13}{5}$. \hfill $\Box$\\

Due to Lemma \ref{exponormperfeisen}, we need to apply Lemma \ref{Ineqodddiv} to only five residue classes. Fortunately, some of them can be dealt with at the same time. The following lemma will now join most of the results we have acquired so far.

\begin{lemma}\label{kequivpm1}
Let $\alpha = (1 - \omega^2)^{k-1}\mu \in \mathbb{Z}[\omega]$ be norm-perfect with $\mu$ odd and a rational $k \geq 2$. Then $k \equiv \pm 1 \mod 12\mathbb{Z}$, $(1 - \omega^2)^k - 1$ is prime and either $(1 - \omega^2)^k - 1$ or $\overline{(1 - \omega^2)^k - 1}$ divide $\mu$.
\end{lemma}

\textit{Proof.} Let $\psi$ be an odd positive prime divisor of
$$(1 - \omega^2)^k - 1 \simeq \sigma_{\mathbb{Q}(\omega)}\left((1 - \omega^2)^{k-1}\right),$$
say $(1 - \omega^2)^k - 1 = \psi\rho$ for a $\rho \in \mathbb{Z}[\omega]$, and such that $\psi$ has the least norm among such divisors of $(1 - \omega^2)^k - 1$. Since $\alpha$ is norm-perfect and the norm of any unit is 1, we have
$$3N_{\mathbb{Q}(\omega)}(\alpha) = N_{\mathbb{Q}(\omega)}\left(\sigma_{\mathbb{Q}(\omega)}(\alpha)\right) = N_{\mathbb{Q}(\omega)}\left(\sigma_{\mathbb{Q}(\omega)}(\mu)\right)N_{\mathbb{Q}(\omega)}(\rho)N_{\mathbb{Q}(\omega)}(\psi).$$
Since $\psi$ is an odd prime, it follows either $\psi \, | \, \alpha$ or $\overline{\psi} \, | \, \alpha$. It is $N_{\mathbb{Q}(\omega)}(\psi) = N_{\mathbb{Q}(\omega)}(\overline{\psi})$, so we may apply the Lemma \ref{Ineqodddiv}. We examine the norm of $(1 - \omega^2)^k - 1$ depending on $k$.
\begin{table}[H]
    \centering
    $$\begin{array}{c|c|c}
        k \mod 12\mathbb{Z} & \Re\left((1 - \omega^2)^k\right) & N_{\mathbb{Q}(\omega)}\left((1 - \omega^2)^k - 1\right) \\
        \hline
        \pm 2 & \frac{3^{\frac{k}{2}}}{2}& 3^k - 3^{\frac{k}{2}} + 1\\
        \pm 1 & \frac{3^{\frac{k + 1}{2}}}{2}& 3^k - 3^{\frac{k + 1}{2}} + 1 \\
        0 & 3^{\frac{k}{2}}& 3^k - 2 \cdot 3^{\frac{k}{2}} + 1
    \end{array}$$
    \caption{The norm of $(1 - \omega^2)^k - 1$ depending on $k$.}
    \label{Normof(1-omega)^k}
\end{table}
Notice the symmetry in this table. We go ahead and insert these numbers into our inequality from Lemma \ref{Ineqodddiv}.
\begin{table}[H]
    \centering
    $$\begin{array}{c|c}
        k \mod 12\mathbb{Z} & N_{\mathbb{Q}(\omega)}(\psi) > \\
        \hline
        \pm 2 & \frac{13}{5} \cdot \frac{3^k - 3^{\frac{k}{2}} + 1}{3^{\frac{k}{2}} - 1}\\
        \pm 1 & \frac{13}{5} \cdot \frac{3^k - 3^{\frac{k + 1}{2}} + 1}{3^{\frac{k + 1}{2}} - 1}\\
        0 & \frac{13}{5} \cdot \frac{3^k - 2 \cdot 3^{\frac{k}{2}} + 1}{2 \cdot 3^{\frac{k}{2}} - 1}
    \end{array}$$
    \caption{Some estimates for $N_{\mathbb{Q}(\omega)}(\psi)$}
    \label{EstimatesForNormPsi}
\end{table}

We will distinguish between three cases.

\underline{\textit{Case 1:}} Suppose $k \equiv 0 \mod 12\mathbb{Z}$, say $k = 12t$ for some $t \in \mathbb{N}_{>0}$. Then, the norm of $(1 - \omega^2)^k - 1$ is
$$(3^{6t} - 1)^2 = (3^t - 1)^2(3^t + 1)^2(3^{2t} - 3^t + 1)^2(3^{2t} + 3^t + 1)^2.$$
As $\psi$ is a prime divisor of $\alpha$ of minimal norm, we have $N_{\mathbb{Q}(\omega)}(\psi) \leq (3^t - 1)^2$. However, Table \ref{EstimatesForNormPsi} gives
\begin{align*}
    N_{\mathbb{Q}(\omega)}(\psi) &> \frac{13}{5} \cdot \frac{3^k - 2 \cdot 3^{\frac{k}{2}} + 1}{2 \cdot 3^{\frac{k}{2}} - 1}\\
    &= \frac{13}{5} \cdot \frac{(3^{\frac{k}{2}} - 1)^2}{2 \cdot 3^{\frac{k}{2}} - 1}\\
    &= \frac{13}{10} \cdot \frac{(3^{\frac{k}{2}} - 1)^2}{3^{\frac{k}{2}} - \frac{1}{2}}\\
    &> 3^{\frac{k}{2}} - 1\\
    &> (3^t - 1)^2
\end{align*}
for $t \geq 1$.

\underline{\textit{Case 2:}} For $k \equiv \pm 2 \mod 12\mathbb{Z}$, we have
$$N_{\mathbb{Q}(\omega)}(\psi) > \frac{13}{5} \cdot \frac{3^k - 3^{\frac{k}{2}} + 1}{3^{\frac{k}{2}} - 1} > \frac{13}{5} \cdot \frac{3^k - 2 \cdot 3^{\frac{k}{2}} + 1}{3^{\frac{k}{2}} - 1} = \frac{13}{5} (3^{\frac{k}{2}} - 1).$$
We will show a more direct contradiction for both cases.

\textit{Subcase 2.1:} Suppose $k \equiv 2 \mod 12\mathbb{Z}$, say $k = 12t + 2$ for some $t \in \mathbb{N}_{\geq0}$. Then, the norm of $(1 - \omega^2)^k - 1$ is
$$3^{12t + 2} - 3^{6t + 1} + 1 = (3^{6t + 1} - 3^{3t + 1} + 1)(3^{6t + 1} + 3^{3t + 1} + 1).$$
The left factor is less than $\frac{13}{5} (3^{\frac{k}{2}} - 1)$ since $\frac{k}{2} = 6t + 1$.

\textit{Subcase 2.2:} Suppose $k \equiv -2 \mod 12\mathbb{Z}$, say $k = 12t - 2$ for some $t \in \mathbb{N}_{>0}$. Then, the norm of $(1 - \omega^2)^k - 1$ is
$$3^{12t - 2} - 3^{6t - 1} + 1 = (3^{6t - 1} - 3^{3t} + 1)(3^{6t - 1} + 3^{3t} + 1).$$
The left factor is less than $\frac{13}{5} (3^{\frac{k}{2}} - 1)$ since $\frac{k}{2} = 6t - 1$.

\underline{\textit{Case 3:}} Suppose $k \equiv \pm 1 \mod 12\mathbb{Z}$ such that $k \geq 11$. By Table \ref{EstimatesForNormPsi}, we get
\begin{align*}
    N_{\mathbb{Q}(\omega)}(\psi) &> \frac{13}{5} \cdot \frac{3^k - 3^{\frac{k + 1}{2}} + 1}{3^{\frac{k + 1}{2}} - 1}\\
    &= \frac{13}{5} \cdot \frac{3^{k-1} - 3^{\frac{k - 1}{2}} + \frac{1}{3}}{3^{\frac{k - 1}{2}} - \frac{1}{3}}\\
    &> \frac{13}{5} \cdot \frac{3^{k-1} - 2 \cdot 3^{\frac{k - 1}{2}} + 1}{3^{\frac{k - 1}{2}} - \frac{1}{3}}\\
    &> \frac{12}{5} \cdot \frac{3^{k-1} - 2 \cdot 3^{\frac{k - 1}{2}} + 1}{3^{\frac{k - 1}{2}} - 1}\\
    &= \frac{12}{5} \left(3^{\frac{k-1}{2}} - 1\right)\\
    &> \sqrt{N_{\mathbb{Q}(\omega)}\left((1 - \omega^2)^k - 1\right)},
\end{align*}
so, by multiplicity of the norm and minimality of $N_{\mathbb{Q}(\omega)}(\psi)$, $(1 - \omega^2)^k - 1$ has exactly one prime divisor, hence it is prime and $\psi = (1 - \omega^2)^k - 1$ or $\psi = \overline{(1 - \omega^2)^k - 1}$. Thus, either $(1 - \omega^2)^k - 1$ or $\overline{(1 - \omega^2)^k - 1}$ divide $\mu$.

The last inequality in the previous computation holds because
\begin{align*}
    \left(\frac{12}{5}\left(3^{\frac{k-1}{2}} - 1\right)\right)^2 &> \left(\frac{11}{5}\left(3^{\frac{k-1}{2}} - \frac{1}{2}\right)\right)^2\\
    &= \frac{121}{25}\left(3^{k-1} - 3^{\frac{k-1}{2}} + \frac{1}{4}\right)\\
    &= \frac{121}{75}\left(3^k - 3^{\frac{k+1}{2}} + \frac{3}{4}\right)\\
    &> 3^k - 3^{\frac{k+1}{2}} + 1\\
    &= N_{\mathbb{Q}(\omega)}((1 - \omega^2)^k - 1)
\end{align*}
for $k \geq 11$. This finishes the proof. \hfill$\Box$\\

We are ready to take the last few steps towards proving the main theorem of this note.

\begin{lemma}
Let $\alpha = (1 - \omega^2)^{k-1}\mu$ be a norm-perfect \textsc{Eisenstein} integer such that $k \geq 2$ and $\mu$ is odd.
\begin{enumerate}
    \item If $k \equiv 1 \mod 12\mathbb{Z}$, then $\mu = \left((1 - \omega^2)^k - 1\right)\rho$ for some odd $\rho$, and
    \item if $k \equiv -1 \mod 12\mathbb{Z}$, then $\mu = \left(\overline{(1 - \omega^2)^k - 1}\right)\rho$ for some odd $\rho$.
\end{enumerate}
\end{lemma}

\textit{Proof.} Suppose $k \equiv 1 \mod 12\mathbb{Z}$. Since $k \geq 2$, $\overline{(1 - \omega^2)^k - 1}$ is not a unit and its positive associate is
$$(1 + \omega)(\overline{(1 - \omega^2)^k - 1}) = 3(1 - \omega)^{k-2} - 1 - \omega.$$
Thus,
$$\sigma_{\mathbb{Q}(\omega)}(\overline{(1 - \omega^2)^k - 1}) = 3(1 - \omega)^{k-2} - \omega$$
because $(1 - \omega^2)^k - 1$ and its complex conjugate are prime by Lemma \ref{kequivpm1}. Furthermore,
\begin{align*}
    N_{\mathbb{Q}(\omega)}\left(\sigma_{\mathbb{Q}(\omega)}(\overline{(1 - \omega^2)^k - 1})\right) &= \left(3(1 - \omega)^{k-2} - \omega\right)\left(3(1 - \omega^2)^{k-2} - \omega^2\right)\\
    &= 3^k - 6\Re(\omega(1 - \omega^2)^{k-2}) + 1.
\end{align*}
We compute, since $k$ is a \textit{rationally} odd natural number, that
$$\omega(1 - \omega^2)^{k-2} = (\omega - 1)^{k-2} = -\left(3(1 - \omega)^{k-3}\right),$$
so
$$\Re(\omega(1 - \omega^2)^{k-2}) = - \frac{3^{\frac{k - 1}{2}}}{2}.$$
Thus,
$$N_{\mathbb{Q}(\omega)}\left(\sigma_{\mathbb{Q}(\omega)}(\overline{(1 - \omega^2)^k - 1})\right) = 3^k + 3^{\frac{k+1}{2}} + 1 > 3^k = 3N_{\mathbb{Q}(\omega)}\left((1 - \omega^2)^{k-1}\right).$$
Since
$$N_{\mathbb{Q}(\omega)}\left(\sigma_{\mathbb{Q}(\omega)}((1 - \omega^2)^{k-1})\right) = N_{\mathbb{Q}(\omega)}\left((1 - \omega^2)^k - 1\right) = N_{\mathbb{Q}(\omega)}\left(\overline{(1 - \omega^2)^k - 1}\right)$$
their product is greater than $3N_{\mathbb{Q}(\omega)}\left((1 - \omega^2)^{k-1}(\overline{(1 - \omega^2)^k - 1})\right)$, hence the prime dividing $\mu$ must be $(1 - \omega^2)^k - 1$.

The proof for $k \equiv - 1 \mod 12\mathbb{Z}$ works analogously.\hfill $\Box$\\

This lemma proves in particular that a norm-perfect number is never divisible by $3^k - 3^{\frac{k + 1}{2}} + 1$, i.e. by the product of an \textsc{Eisenstein Mersenne} prime and its conjugate. The rest is just putting all the pieces together.\\

\textit{Proof of Main Theorem \ref{EuclidEulerEisenstein}}. We start with the proof about the even norm-perfect numbers. The \textit{only-if}-direction is a corollary of the previous lemma since we proved that every norm-perfect number is divisible by exactly one of the presented primitive norm-perfect numbers.

The \textit{if}-direction is a simple computation and will be shown for the case $k \equiv -1 \mod 12\mathbb{Z}$. Quickly check beforehand that $\overline{(1 - \omega^2)^k - 1}$ is positive. It is
\begin{align*}
    N_{\mathbb{Q}(\omega)}\left(\sigma_{\mathbb{Q}(\omega)}(\alpha)\right) &= N_{\mathbb{Q}(\omega)}\left(\frac{(1 - \omega^2)^k - 1}{1 - \omega^2 - 1}\cdot \left( \left(\overline{(1 - \omega^2)^k - 1}\right) + 1\right)\right)\\
    &= N_{\mathbb{Q}(\omega)}\left((1 - \omega^2)^k - 1\right)N_{\mathbb{Q}(\omega)}\left( \left(\overline{(1 - \omega^2)^k - 1}\right) + 1\right)\\
    &= N_{\mathbb{Q}(\omega)}\left((1 - \omega^2)^k - 1\right)N_{\mathbb{Q}(\omega)}\left(\overline{(1 - \omega^2)^k}\right)\\
    &= 3N_{\mathbb{Q}(\omega)}\left((1 - \omega^2)^k - 1\right)N_{\mathbb{Q}(\omega)}\left((1 - \omega^2)^{k-1}\right)\\
    &= 3N_{\mathbb{Q}(\omega)}\left(\alpha\right),
\end{align*}
hence $\alpha$ is norm-perfect. However, $\overline{(1 - \omega^2)^k - 1}$ and any power of $(1 - \omega^2)$ are deficient, so $\alpha$ is primitive.

The computation for $k \equiv 1 \mod 12\mathbb{Z}$ works analogously.

Now we prove the case of even perfect numbers. The \textit{if}-direction is given by a similar computation as above but without the norm function. Notice, if $k \equiv 1 \mod 12\mathbb{Z}$, then
$$\sigma_{\mathbb{Q}(\omega)}\left((1 - \omega^2)^k - 1\right) = (1 - \omega^2)^k,$$
yielding
$$\sigma_{\mathbb{Q}(\omega)}\left((1 - \omega^2)^{k-1}((1 - \omega^2)^k - 1)\right) = (1 - \omega^2)^k((1 - \omega^2)^k - 1).$$

Conversely, the \textit{only-if}-direction is proven by examining all the possible primitive norm-perfect numbers, since every perfect number is also norm-perfect. The associates of a perfect number cannot be perfect, since all of them have the same image under the $\sigma_{\mathbb{Q}(\omega)}$-function. In the case of $k \equiv -1 \mod 12\mathbb{Z}$, we see from the computations in the previous proof that the equation
$$\sigma_{\mathbb{Q}(\omega)}(\alpha) = (1 - \omega^2)\alpha$$
cannot hold, since the left-hand side is divisible by $(1 - \omega^2)^k - 1$ and the right-hand side by its complex conjugate but neither of them is divisible by the product of those two factors, which is $3^k - 3^{\frac{k + 1}{2}} + 1$.

Finally, by Lemma \ref{expomersennecyclo}, we know that the primality of $(1 - \omega)^k - 1$ implies the same for $k$.
\hfill $\Box$

\subsection{Odd \textsc{Eisenstein} integers}

Similarly to the previous cases of the rational and the \textsc{Gaussian} integers, the case of odd norm-perfect \textsc{Eisenstein} integers proves to be a bit more difficult to come by. \textsc{Parker} et al. presented a form for such a number in their conjecture \ref{conPark}. We will work in a similar way like \textsc{Ward} who transferred Theorem \ref{Euler} to $\mathbb{Z}[i]$. Due to two congruence classes of odd numbers being existent in $\mathbb{Z}[\omega]$, we will have to start a bit differently.

\begin{lemma}
\label{oddeisen}
Let $\psi \in \mathbb{P}_{\mathbb{Q}(\omega)}^+$ be an odd positive prime and $m \in \mathbb{N}$. It is $$N_{\mathbb{Q}(\omega)}(\sigma_{\mathbb{Q}(\omega)}(\psi^m)) \equiv 0 \mod 3\mathbb{Z}$$
if and only if one of the following cases applies:
\begin{enumerate}
    \item $\psi \equiv 1 \mod 1 - \omega^2$ and $m \equiv 2 \mod 3\mathbb{Z}$.
    \item $\psi \equiv 2 \mod 1 - \omega^2$ and $m \equiv 1 \mod 2\mathbb{Z}$.
\end{enumerate}
\end{lemma}

\textit{Proof.} As 3 totally ramifies in $\mathbb{Q}(\omega)$, the norm of $\sigma_{\mathbb{Q}(\omega)}(\psi^m) = \sum_{k = 0}^m \psi^k$ is divisible by 3 if and only if the sum itself is even. Recall that
$$\nicefrac{\mathbb{Z}[\omega]}{\langle 1 - \omega^2\rangle} \cong \nicefrac{\mathbb{Z}}{3\mathbb{Z}}.$$
\begin{enumerate}
    \item If $\psi \equiv 1 \mod 1 - \omega^2$, each of its powers is congruent to $1 \mod 1 - \omega^2$. Thus, 
    $$\sum_{k = 0}^m \psi^k \equiv m + 1 \mod 1 - \omega^2$$
    which is 0 if and only if $m \equiv 2 \mod 3\mathbb{Z}$.
    \item If $\psi \equiv 2 \mod 1 - \omega^2$, we have $\psi^l \equiv l \mod 2\mathbb{Z}$. For \textit{rationally} odd $k \in \mathbb{N}$, this implies
    $$\psi^{k-1} + \psi^k \equiv 0 \mod 1 - \omega^2.$$
    Now just use the sum above. \hfill $\Box$\\
\end{enumerate}

After comparing this result to the case of the rational or \textsc{Gaussian} integers, we see that we the number of possible forms increases. Having this lemma in our repertoire, we may now also prove Theorem \ref{oddnormperfEisenstein}.\\

\textit{Proof of Theorem \ref{oddnormperfEisenstein}.} Since $\alpha$ is norm-perfect, we have
$$N_{\mathbb{Q}(\omega)}\left(\sigma_{\mathbb{Q}(\omega)}(\alpha)\right) = 3N_{\mathbb{Q}(\omega)}(\alpha).$$
Write
$$\alpha = \varepsilon\prod_{\psi_1 \in P_1}\psi_1^{e_{\psi_1}}\prod_{\psi_2 \in P_2}\psi_2^{e_{\psi_2}}.$$
As the norm and $\sigma_{\mathbb{Q}(\omega)}$-function are both multiplicative, so is their composition and we may write the first equation as
$$\prod_{\psi_1 \in P_1}N_{\mathbb{Q}(\omega)}\left(\sigma_{\mathbb{Q}(\omega)}(\psi_1^{e_{\psi_1}})\right)\prod_{\psi_2 \in P_2}N_{\mathbb{Q}(\omega)}\left(\sigma_{\mathbb{Q}(\omega)}(\psi_2^{e_{\psi_2}})\right)$$
$$= 3\prod_{\psi_1 \in P_1}N_{\mathbb{Q}(\omega)}\left(\psi_1^{e_{\psi_1}}\right)\prod_{\psi_2 \in P_2}N_{\mathbb{Q}(\omega)}\left(\psi_2^{e_{\psi_2}}\right),$$
using that the norm of a unit is 1 and the $\sigma_{\mathbb{Q}(\omega)}$-function is defined up to associates. By Lemma \ref{lemminimalprime}, the norm of any odd \textsc{Eisenstein} is congruent to 1 modulo $3\mathbb{Z}$. Therefore, the right-hand side is congruent to 3 modulo $9\mathbb{Z}$ implying that exactly one of the prime factors of $\alpha$ satisfies the conditions listed in Lemma \ref{oddeisen}. The possible forms are the ones presented above. \hfill $\Box$\\

Theorem \ref{oddnormperfEisenstein} is a bit weaker than \textsc{Parker} et al.'s conjecture because the latter corresponds to the case of the former where $\psi_0 \in P_1$. Nevertheless, we get a corollary if we distinguish between the two congruence classes of odd integers.

\begin{corollary}
In the case of Theorem \ref{oddnormperfEisenstein}, each case corresponds to exactly one congruence class of odd positive integers, i.e. $\psi_0 \in P_1$ if and only if $\alpha \equiv 1 \mod 1 - \omega^2$ under the assumption that $\alpha$ is positive.
\end{corollary}

\textit{Proof.} We evaluate the product from the theorem and remark that $\varepsilon = 1$ since we are focusing on positive integers. In any case, the powers
$$\psi_1^{e_{\psi_1}} \equiv 1 \mod 1 - \omega^2$$
since $\psi_1 \equiv 1 \mod 1 - \omega^2$ and also
$$\psi_2^{e_{\psi_2}} \equiv 1 \mod 1 - \omega^2$$
since $e_{\psi_2} \equiv 0 \mod 2\mathbb{Z}$ and $\psi_2 \equiv 2 \mod 1 - \omega^2$. Hence, the congruence class of the whole product only depends on $\psi_0^k$ which is $\equiv 1$ modulo $1 - \omega^2$ if $\psi_0 \in P_1$ and $\equiv 2$ modulo $1 - \omega^2$ if $\psi_0 \in P_2$.\hfill $\Box$\\

Similarly to the case of the rational integers, neither the author nor \textsc{Parker} et al. \cite{Parker} have been able to present any odd norm-perfect integers yet.

\begin{remark}
One might hope to find easy solutions when considering the rational perfect integers because $\mathbb{Z} \subset \mathbb{Z}[\omega]$ and 2 is an odd \textsc{Eisenstein} prime. Let $N = 2^{k-1}(2^k - 1)$ be a perfect rational integer, so $k$ is a rational prime and, in order for $N$ to be odd in the \textsc{Eisenstein} integers, we assume that $k > 2$. Theorem \ref{oddnormperfEisenstein} implies that $\psi_0$ must be an \textsc{Eisenstein} prime above $2^k - 1$ because $k$ is rationally odd. However, $2^k - 1$ splits in $\mathbb{Q}(\omega)$ because $2^k - 1 \equiv 1 \mod 3\mathbb{Z}$ due to $k \equiv 1 \mod 2 \mathbb{Z}$ and thus the two \textsc{Eisenstein} primes above it are $\psi_0$ and $\overline{\psi_0}$ which also lie in the same residue class modulo $1 - \omega^2$. Since both also have the same exponent---namely 1---, this yields a contradiction and $N$ cannot be norm-perfect.
\end{remark}

Nevertheless, we can prove that the \textsc{Gaussian} integers feature a kind of numbers which cannot be found in the \textsc{Eisenstein} integers.

\begin{theorem}
 There is no norm-perfect \textsc{Eisenstein} prime.
\end{theorem}

\textit{Proof.}
Recall that $1 - \omega^2$ is not norm-perfect by Main Theorem \ref{EuclidEulerEisenstein}. Furthermore, suppose $\psi$ is an odd norm-perfect positive \textsc{Eisenstein} prime. Then,
$$3N_{\mathbb{Q}(\omega)}(\psi) = N_{\mathbb{Q}(\omega)}\left(\sigma_{\mathbb{Q}(\omega)}(\psi)\right) = N_{\mathbb{Q}(\omega)}(\psi + 1).$$
Evaluating this yields the equation
$$3(a^2 - ab + b^2) = (a + 1)^2 - (a + 1)b + b^2$$
which has the integer solutions $a = 0, b = - 1$ and $a = b = 1$ but $-\omega$ and $1 + \omega$ are units and non-positive, too. Since all associates of a norm-perfect number are also norm-perfect, this is sufficient to prove the claim. \hfill $\Box$\\

An interesting fact is that the previous theorem is independent of our choice of the set of positive primes.

The results in this subsection apply solely to odd norm-perfect integers but the author believes that there is not much more that can be said about odd perfect integers and has not been listed here.

\subsection{Cyclotomic fields of higher degree}

Recall our set $\mathcal{R}$ of cyclotomic fields over $\mathbb{Q}$ with class number 1 and generated by $\zeta_p$ with $p$ being a rational prime or 4. So far, we have presented results concerning the two fields with the smallest degree, $\mathbb{Q}(i)$ and $\mathbb{Q}(\omega)$. We will have a look at the other fields in $\mathcal{R}$ and there is a nice theorem about which fields exactly are contained in it.

\begin{theorem}[{{\cite[theorem 11.1]{Washington}}}] The numbers $n \in \mathbb{N}$ such that $\mathbb{Q}(\zeta_n)$ has class number 1 are:
\begin{itemize}
    \item 1 through 22, 24, 25, 26, 27, 28, 30, 32, 33, 34, 35, 36, 38, 40, 42, 44, 45, 48, 50, 54, 60, 66, 70, 84, and 90.
\end{itemize}
\end{theorem}

The set $\mathcal{R}$ does not contain all of those fields but the number being finite tells us that the problem of (norm-)perfect numbers in such fields is a soluble one. It is even less than we may think in the first moment, since it is $\mathbb{Q}(\zeta_n) = \mathbb{Q}(\zeta_{2n})$ for odd $n \in \mathbb{N}$. Thus, our shortened list of $n$ such that $\mathbb{Q}(\zeta_n)$ is a UFD is:
\begin{itemize}
    \item 2, 3, 4, 5, 7, 8, 9, 11, 12, 13, 15, 16, 17, 19, 20, 21, 24, 25, 27, 28, 32, 33, 35, 36, 40, 44, 45, 48, 60, and 84.
\end{itemize}
Picking $p$ such that $p$ is a rational prime or 4, yields the list
\begin{itemize}
    \item 2, 3, 4, 5, 7, 11, 13, 17, and 19.
\end{itemize}
The cases $p \in \{2,3,4\}$ were dealt with in the introduction and preceding subsections. Since the degree $[\mathbb{Q}(\zeta_p):\mathbb{Q}] = \varphi(p)$, the norm of an element in $\mathbb{Z}[\zeta_p]$ is the product of an increasing number of factors which complicates its computation. However, based on our experiences already gained, we may formulate a conjecture.

\begin{conjecture}
Let $K = \mathbb{Q}(\zeta_p) \in \mathcal{R}$ with $p \neq 2, 4$ where $\zeta_p$ is the primitive $p$-th root of unity such that $1 - \zeta_p$ is positive in $\mathbb{Z}[\zeta_p]$. Then,
\begin{itemize}
    \item the even perfect integers in this ring are given by
    $$-\zeta_p^{-1}(1 - \zeta_p)^{k-1}\left((1 - \zeta_p)^k - 1\right)$$
    for $k \equiv 1 \mod 4p\mathbb{Z}$ and $(1 - \zeta_p)^k - 1$ being prime and
    \item the even norm-perfect integers are the associates of the perfect integers from item 1 and the associates of
    $$(1 - \zeta_p)^{k-1}\left(\overline{(1 - \zeta_p)^k - 1}\right)$$
    for $k \equiv -1 \mod 4p\mathbb{Z}$ and $(1 - \zeta_p)^k - 1$ being prime.
\end{itemize}
\end{conjecture}

When we turn to the odd norm-perfect integers in $\mathbb{Z}[\zeta_p]$, however, presenting the form they have is possible. As
$$\nicefrac{\mathbb{Z}[\zeta_p]}{\langle 1 - \zeta_p \rangle} \cong \nicefrac{\mathbb{Z}}{p\mathbb{Z}},$$
we have to take care about an increasing amount of residue classes of odd primes which may be factors of an odd norm-perfect number $\alpha$. Thus, the simplicity of the form of $\alpha$ heavily depends on the structure of $\left(\nicefrac{\mathbb{Z}}{p\mathbb{Z}}\right)^*$.

\begin{lemma}
Let $p$ be an odd natural prime, $a \in \mathbb{N}$ such that $a \not \equiv 0, 1 \mod p\mathbb{Z}$ and $t$ the order of $a$ in $\left(\nicefrac{\mathbb{Z}}{p\mathbb{Z}}\right)^*$. Then $p \, | \, \sum_{k = 0}^{t - 1} a^k$.
\end{lemma}

\textit{Proof.} By assumption $t > 1$. Then $a$ is a root of
$$x^t - 1 = (x - 1)\sum_{k = 0}^{t - 1}x^k \mod p\mathbb{Z}.$$
Since $a \not \equiv 1 \mod p\mathbb{Z}$ and $\nicefrac{\mathbb{Z}}{p\mathbb{Z}}$ is a field, $a$ has to be a root of the sum, and we get $\sum_{k = 0}^{t - 1} a^k \equiv 0 \mod p\mathbb{Z}$. \hfill $\Box$\\

The following theorem is the generalisation of the work we got familiar with while working with odd norm-perfect \textsc{Eisenstein} integers.

\begin{theorem}
Let $K = \mathbb{Q}(\zeta_p) \in \mathcal{R}$ with $p \neq 2, 4$ where $\zeta_p$ is the primitive $p$-th root of unity such that $1 - \zeta_p$ is positive in $\mathbb{Z}[\zeta_p]$. Let $\alpha \in \mathbb{Z}[\zeta_p]$ be an odd norm-perfect integer and
$$P_j = \{\psi \in \mathbb{P}_{\mathbb{Q}(\zeta_p)}^+: \psi \equiv j \mod 1 - \zeta_p \land \psi \, | \, \alpha\}$$
for all $j \in \{1, \ldots, p - 1\}$. Then $\alpha$ is of the form
$$\alpha = \varepsilon \psi_0^k \prod_{j = 1}^{p - 1}\left(\prod_{\psi_j \in P_j \setminus \{\psi_0\}} \psi_j^{e_{\psi_j}}\right)$$
where
\begin{enumerate}
    \item each $\psi$ with a subscript is an odd positive prime,
    \item $\varepsilon$ is a unit,
    \item if $j \neq 1$, $e_{\psi_j} \not \equiv -1 \mod t\mathbb{Z}$ where $t$ is the order of $j$ in $\left(\nicefrac{\mathbb{Z}}{p\mathbb{Z}}\right)^*$,
    \item $e_{\psi_1} \not \equiv -1 \mod p\mathbb{Z}$, and
    \item if $\psi_0 \in P_j$ for $j \neq 1$, then $k \equiv -1 \mod t\mathbb{Z}$ where $t$ is the order of $j$ in $\left(\nicefrac{\mathbb{Z}}{p\mathbb{Z}}\right)^*$, or if $\psi_0 \in P_1$, then $k \equiv -1 \mod p\mathbb{Z}$.
\end{enumerate}
\end{theorem}

\textit{Proof.} We use that, since $\alpha$ is odd, $N_K(\alpha)$ is not divisible by $p$ because $p$ is totally ramified in $\mathbb{Q}(\zeta_p)$. Thus, $p^2 \, \nmid \, p N_K(\alpha) = N_K\left(\sigma_K(\alpha)\right)$. By the same argument as in the proof of Theorem \ref{oddnormperfEisenstein}, there is only one particular positive prime divisor $\psi_0$ of $\alpha$ such that $p \, | \, N_K\left(\sigma_K(\psi_0^k)\right)$ for $k$ being the exponent of $\psi_0$ in the prime decomposition of $\alpha$. For $\psi_0 \in P_1$ this is the case if and only if $k \equiv -1 \mod p\mathbb{Z}$ and, for $j \neq 1$ and $\psi_0 \in P_j$, this is the case if and only if $k \equiv -1 \mod t\mathbb{Z}$ by the previous lemma. All the other primes must not satisfy these congruence conditions. \hfill $\Box$


\subsection*{Acknowledgements}
My thanks go to Victoria Cantoral Farfán who supervised my master's thesis, encouraged me to submit its first part as this note, and guided me in the process.


\end{document}